\documentclass[12pt,reqno,intlim]{amsart}

\usepackage{amssymb,amsmath}

\usepackage[arrow,matrix,curve]{xy}
\newdir{ >}{{}*!/-5pt/\dir{>}}

\newcommand{\nc}{\newcommand}

\nc{\bC}{\bold{C}} \nc{\bN}{\Bbb{N}} \nc{\cF}{\mathcal{F}}
\nc{\cE}{\mathcal{E}} \nc{\cR}{\mathcal{R}} \nc{\cM}{\mathcal{M}}
\nc{\al}{\alpha} \nc{\bt}{\beta} \nc{\gm}{\gamma} \nc{\dl}{\delta}
\nc{\om}{\omega} \nc{\sg}{\sigma} \nc{\Sg}{\Sigma} \nc{\vf}{\varphi}
\nc{\ve}{\varepsilon} \nc{\os}{\overset} \nc{\ol}{\overline}
\nc{\ul}{\underline} \nc{\us}{\underset} \nc{\sbs}{\subset}
\nc{\bsl}{\backslash} \nc{\Ra}{\Rightarrow}
\nc{\lra}{\longrightarrow} \nc{\all}{\allowdisplaybreaks}

\nc{\Codes}{\operatorname{{\bold{Codes}}}}
\nc{\RegMono}{\operatorname{\mathcal{R}{\rm{eg}\mathcal{M}{\rm{ono}\!}}}}
\nc{\RegEpi}{\operatorname{\mathcal{R}{\rm{eg}\mathcal{E}{\rm{pi}\!}}}}
\nc{\Mn}{\operatorname{\mathcal{M}{\rm{ono}\!}}}
\nc{\Ep}{\operatorname{\mathcal{E}{\rm{pi}\!}}}
\nc{\Rg}{\operatorname{\mathcal{R}{\rm{eg}\!}}}
\nc{\Ob}{\operatorname{Ob\!}}

\numberwithin{equation}{section}

\newtheorem{theo}{\ \ \ Theorem}[section]
\newtheorem{lem}[theo]{\ \ \ Lemma}
\newtheorem{prop}[theo]{\ \ \ Proposition}
\newtheorem{cor}[theo]{\ \ \ Corollary}

\theoremstyle{definition}
\newtheorem{exmp}[theo]{\ \ \ Example}

\theoremstyle{remark}
\newtheorem{rem}[theo]{\ \ \ Remark}

\begin{document}

\title[]
{The factorization system of a radical on a homological category}

\author{Dali Zangurashvili}

\maketitle

% Abstract.
\begin{abstract}
Employing the techniques of transporting a factorization system from one category to another  via an adjunction, developed in our earlier papers, a factorization system is related to a radical on a homological category. 
%A positive answer is given in two cases:  (a) where the ground category is complete, well-powered, and homological, and a radical is arbitrary, (b) where the ground category is an arbitrary homological one, but a radical is idempotent.

\vskip+2mm
\noindent{\bf Key words and phrases}: factorization system, homological category, radical, idempotent radical, semi-left-exact reflection.
\vskip+2mm

\end{abstract}

% 1.
\section{Introduction}
In view of the well-known theorem on reflective factorization systems proposed by Cassidy, H\'{e}bert, and Kelly \cite{CHK}, the problem of transporting a factorization system from one category to another  via an adjunction was studied in \cite{Z1}-\cite{Z3}. 
In the present paper, we employ the techniques developed in these papers to construct a  factorization system  on a homological category provided that a radical on it is given. We use the term `homological category' in the sense of the book \cite{BB} by Borceux and Bourn, i.e., for a pointed regular protomodular category. The term was suggested by Johnstone, as is noted by the authors, and provided by the fact that the classical homological lemmas (the short five lemma, the nine lemma, the snake lemma), as well as the Noether isomorphism theorems hold true in such categories. The class of homological categories contains that of abelian categories, but also many other interesting categories, including semi-abelian categories introduced by  Janelidze, M\'{a}rki and Tholen \cite{JMT}, such as the categories of groups (more generally, those of $\Omega$-groups, including the categories of associative rings and crossed modules), Heyting semilattices,  compact Hausdorff groups (more generally, the categories of compact Hausdorff $\Omega$-groups), etc. Other examples of homological categories are provided by some quasi-varieties, and also by the categories of topological groups, topological associative rings, and others.

 The notion of a radical on a homological category, used in this paper, is a straightforward  generalization of that of a radical on the category of left/right modules over an associative ring with identity in the sense of classical torsion theory (see, e.g., \cite{L}).

To accomplish our goal, we consider the full refletive subcategory $\mathcal{X}$ of $\mathbf{R}$-torsionfree objects of a homological category $\mathcal{C}$, where $\mathbf{R}$ is a radical on $\mathcal{C}$. It is shown that the category $\mathcal{X}$ also is homological, and hence has the factorization system $$(Regular \; Epis, Monos).$$ Transporting this factorization system via the reflection, we arrive at the following morphism classes: 

--  the class $\mathbb{E}$ of morphisms $f:C\rightarrow C'$ such that the composition $\eta_{C'} f$ is a regular epimorphism ($\eta_{C'}$ denotes the cokernel of the normal monomorphism $\mathbf{R}(C')\rightarrowtail C'$). 
It is proved that if the image of $f$ is a normal subobject of $C'$, then $f\in \mathbb{E}$ if and only if $$Im\; f+\mathbf{R}(C')=C';$$

-- the class $\mathbb{M}$ of monomorphisms $m:C\rightarrow C'$ which satisfy the  following condition: if $D$ is a subobject of $C'$ such that the composition of the inclusion $C\cap D\rightarrowtail D$ with the morphism $\eta_D$ is a regular epimorphism, then $D\subseteq C$.

\vskip+1mm
With the aid of the results obtained in \cite{Z1}-\cite{Z3}, we show that
 the pair $(\mathbb{E}, \mathbb{M})$ of morphism classes is a factorization system in any of the following two cases:

-- $\mathcal{C}$ is complete and well-powered;

-- $\mathbf{R}$ is idempotent.
\vskip+2mm
Note that, in the second case, the class $\mathbb{M}$ has a simpler description: a morphism $m$ lies in $\mathbb{M}$ if and only if it is a monomorphism and $\mathbf{R}(C')\subseteq C$.

Financial support from  Shota Rustaveli  National Science Foundation of Georgia
(Ref.: FR-24-8249) is gratefully acknowledged.

\section{Preliminaries}
We use the term `factorization system' in the sense of the paper \cite{FK} by Freyd and Kelly. Namely, a factorization
system on a category $\mathcal{C}$ is a pair of morphism classes $(\mathbb{E},\mathbb{M})$ such that (a) the classes
$\mathbb{E}$ and $\mathbb{M}$ are closed under composition with isomorphisms, (b) every morphism $f$ admits an $(\mathbb{E},\mathbb{M})$-factorization, i.e., there are morphisms $e\in \mathbb{E}$ and 
and $m\in \mathbb{M}$  with $f=me$, and (c) the diagonalization condition holds, i.e., for any $e\in \mathbb{E}$ and $m\in \mathbb{M}$, one has $e\downarrow  m$ (or, equivalently, $m\uparrow e$), which means that for any commutative diagram 
\begin{equation}
\xymatrix{A\ar[r]^{e}\ar[d]_{\alpha}&B\ar[d]^{\beta}\\
C\ar[r]^{m}&D}
\end{equation}
there is a unique morphisms $\delta:B\rightarrow C$ with $m\delta=\beta$ and $\delta e=\alpha$.
%The classes $\mathbb{E}$ and $\mathbb{M}$ in a factorization system satisfy certain stability properties. In particular, they are closed under composition.

Before continue, let us agree on the notation: for any morphism class $\mathbb{N}$, we will use the symbol $\mathbb{N}^{\uparrow}$ (resp. $\mathbb{N}^{\downarrow}$) for the class of all morphisms $f$ with $f \downarrow n$  (resp. $f \uparrow n$) for all $n\in \mathbb{N}$.

Recall that if $(\mathbb{E},\mathbb{M})$ is a factorization system, then $\mathbb{E}^{\downarrow}=\mathbb{M}$ and $\mathbb{M}^{\uparrow}=\mathbb{E}$. For other properties of factorization systems, we refer the reader to the book \cite{B2}.

Let $\mathcal{C}$ and $\mathcal{X}$ be categories, and 
\begin{equation} \label{2.2}
\mathbf{I}:\mathcal{C}\rightarrow \mathcal{X}
\end{equation}
be a functor. Let $C$ be an object of $\mathcal{C}$, and
\begin{equation}
\mathbf{I}^C:\mathcal{C}/C\rightarrow \mathcal{C}/I(C)
\end{equation}
be the induced functor between the slice-categories. Let
$(\mathbb{E},\mathbb{M})$ be a factorization system on $\mathcal{X}$.

\begin{theo} \cite[Theorem 3.4]{Z3} \label{theo2.1}
Let, for each object
$C$ of $\mathcal{C}$, the functor $\mathbf{I}^C$
have a right adjoint $\mathbf{H}^C$, and let all $\mathbf{H}^C$ be full. Then the pair of morphism classes 
\begin{equation} \label{2.4}
(\mathbf{I}^{-1}(\mathbb{E}), (\mathbf{I}^{-1}(\mathbb{E}))^{\downarrow}).
\end{equation}
 is a factorization system.

\end{theo}

If $\mathbf{I}$ is a fibration, then, as is well-known, the functors $\mathbf{H}^C$ exist and are full and faithful. Moreover, the morphisms from the class $(\mathbf{I}^{-1}(\mathbb{E}))^{\downarrow}$ have a simple description. 

\begin{cor} \cite[Example 1.3(2)]{PT}, \cite[Remark 3.10]{Z3} \label{theo2.2}
Let $\mathbf{I}$ be a fibration. Then the pair (\ref{2.4}) is a factorization system on $\mathcal{C}$, and the class $(\mathbf{I}^{-1}(\mathcal{C}))^{\downarrow}$ coincides with that of Cartesian morphisms over all morphisms from $\mathbb{M}$.
\end{cor}

Let $\mathcal{C}$ be a category with pullbacks, and let $\mathbf{I}$ be an arbitrary functor (\ref{2.2}). Assume that it has a right adjoint $\mathbf{H}$. Let $\eta$ be the unit of the adjunction
\begin{equation} \label{2.5}
\xymatrix{\mathcal{C}\ar@<0.7ex>[r]^{\mathbf{I}}&\mathcal{X}\ar[l]^{\mathbf{H}}}
\end{equation}
 Recall that, in that case, for any $C$, the functor $\mathbf{I}^{C}$ has a right adjoint $\mathbf{H}^C$; it sends an object $g$ of $\mathcal{X}/\mathbf{I}(C)$ to the pullback of $\mathbf{H}(g)$ along $\eta_C$. According to \cite[Lemma 2.2]{Z1}, one has  
 \begin{equation} \label{2.6}
 (I^{-1}(\mathbb{E}))^{\downarrow}=\mathbf{H}(\mathbb{M})^{\uparrow\downarrow}.
 \end{equation}  
 Therefore, Theorem \ref{theo2.1} implies the following statement.
%induces an %adjunction
%\begin{equation}
%\xymatrix{\mathcal{C}/C\ar@<0.7ex>[r]^{\mathbf{I}^C}& \mathcal{X}/\mathbf{I}(C)\ar[l]^{\mathbf{H}^C}}
%\end{equation}
%between the slice-categories.

\begin{theo} \cite[Therem 3.7]{Z2}, \cite[Theorem 3.6]{Z3} \label{th2.1.1}
Let $\mathcal{C}$ have pullbacks, and (\ref{2.5}) be an adjunction. If $\mathbf{H}^{C}$ is full, for any object $C$, then the pair 
\begin{equation} \label{2.7}
(\mathbf{I}^{-1}(\mathbb{E}), \mathbf{H}(\mathbb{M})^{\uparrow\downarrow})
\end{equation}
is a factorization system on $\mathcal{C}$. 
\end{theo}

Let now (\ref{2.5}) be a reflection. It is said to be \textit{semi-left-exact} if for any pullback
\begin{equation}
\xymatrix{A\ar[r]^{g}\ar[d]&X\ar[d]^{f}\\
C\ar[r]^{\eta_C}&\mathbf{I}(C)}
\end{equation}
where $X$ is an object of $\mathcal{X}$, the morphism $\mathbf{I}(g)$ is an isomorphism \cite{CHK}. The equivalent notion is that of an \textit{admissible} reflection by Janelidze \cite{J}. It is defined as a reflection (\ref{2.5}) such that all functors $\mathbf{H}^C$  are full and faithful. Therefore, Theorem \ref{th2.1.1} implies the following statement.

\begin{cor} \cite[Therem 3.7]{Z2}, \cite[Theorem 2.5]{Z1} \label{th2.1}
Let $\mathcal{C}$ have pullbacks, and (\ref{2.5}) be a semi-left-exact reflection.  Then the pair of morphism classes (\ref{2.7})
%$$(\mathbf{I}^{-1}(\mathbb{E}), (I^{-1}(\mathbb{E}))^{\downarrow})$$
is a factorization system on $\mathcal{C}$. 
\end{cor}

A morphism $f$ of $\mathcal{C}$ is called a \textit{trivial covering} in Janelidze's Galois theory if the square 
\begin{equation}
\xymatrix{C\ar[r]^{\eta_C}\ar[d]_{f}&\mathbf{I}(C)\ar[d]^{\mathbf{I}(f)}\\
C'\ar[r]^{\eta_{C'}}&\mathbf{I}(C')}
\end{equation}
is a pullback \cite{J}, \cite{J1}.

\begin{prop} \cite[Lemma 2.3]{Z1}, \cite[Corollary 2.6]{Z1} \label{prop2.2}
Let $\mathcal{C}$ have pullbacks. The class $\mathbf{H}(\mathbb{M})^{\uparrow\downarrow}$ contains the class of all trivial coverings which are mapped into $\mathbb{M}$  by $\mathbf{I}$.  If (\ref{2.5}) is  semi-left-exact, then these two classes coincide.
\end{prop}

Note that the requirement that (\ref{2.5}) is semi-left-exact cannot be removed from the second claim of Proposition \ref{prop2.2} (even if we restrict our consideration only to trivial factorization systems) as is shown by \cite[Example 4.2]{CHK}. 
%However, for an arbitrary reflection, the class of all trivial coverings which are mapped into $\mathbb{M}$  by $\mathbf{I}$ is contained in the class $\mathbf{H}(\mathbb{M})^{\uparrow\downarrow}$, as it follows from . 

Under certain additional conditions on a category $\mathcal{C}$ and the class $\mathbb{M}$, the requirement that the functors $\mathbf{H}^{C}$ are full can be removed from Theorem 2.1 and Theorem 2.3. In particular, we have the following statement.

\begin{theo}  \cite[Theorem 5.4]{Z3} \label{th2.3}
Let $\mathcal{C}$ be a complete and well-powered category, and (\ref{2.5}) be an adjunction. Let $\mathbb{M}\subseteq Mono\; \mathcal{C}$ (where $Mono\; \mathcal{C}$ denotes the class of monomorphisms of $\mathcal{C}$). Then the pair of morphism classes (\ref{2.7})
%$$(\mathbf{I}^{-1}(\mathbb{E}), \mathbf{H}(\mathbb{M})^{\uparrow\downarrow})$$
is a factorization system on $\mathcal{C}$.
\end{theo}

For characterizing morphisms from the class $(I^{-1}(\mathbb{E}))^{\downarrow}$, the following lemma can be employed.
\begin{lem} \cite[Lemma 5.5]{Z3} \label{lem2.4}
Let $\mathcal{C}$ be a category with pullbacks, and $\mathbb{N}$ be a class of morphisms closed under composition and satisfying the following cancellation property: if $e=fg$ and $e\in \mathbb{N}$, then $f\in \mathbb{N}$. Let $(\mathbb{E}',\mathbb{M}')$ be a factorization system on $\mathcal{C}$ with $\mathbb{E'}\subseteq \mathbb{N}$. Then $\mathbb{N}^{\downarrow}\subseteq \mathbb{M}'$, and a morphism $m$ from the class $\mathbb{M}'$ lies in $\mathbb{N}^{\downarrow}$ if and only if the following condition is satisfied: if the pullback $\alpha$ of $m$ along a morphism from $\mathbb{M}'$ lies in $\mathbb{N}$, then $\alpha$ is an isomorphism.
\end{lem}

%\begin{rem}
%For the case of an arbitrary reflection, the description of the class %$H(\mathbb{M})^{\uparrow\downarrow}$ is provided by \cite[Lemma 5.5]{Z3}. 
%\end{rem}

Let $\mathcal{C}$ be a pointed category (i.e., a category where an initial and a terminal objects exist and are isomorphic; they are denoted by $0$ and are called zero objects). For such a category, and any objects $C$ and $C'$ from it, the composition $C\rightarrow 0\rightarrow C'$ is denoted by $0_{C.C'}$ and is called a zero morphism. If the product $C\times C'$ exists, then we have the morphisms $l_{C}:C\rightarrow C\times C'$, induced by the pair of morphisms $(1_{C},0_{C,C'})$, and $r_{C'}:C'\rightarrow C\times C'$, induced by $(0_{C',C},1_{C'})$. The morphisms $l_{C}$ and $r_{C'}$ are obviously monomorphisms. 

When no confusion might arise, we skip the indexes in the notation $0_{C,C'}$. 

The term `normal monomorphism' is used for a monomorphism that is a kernel.

Recall that a category $\mathcal{C}$ is called \textit{regular} if it is finitely complete, every kernel pair has a coequalizer and the class of regular epimorphisms is stable under pullback. As is well-known, the pair $(RegEpi \; \mathcal{C},Mono \; \mathcal{C})$ is a factorization system, for any regular category; here $RegEpi \; \mathcal{C}$ denotes the class of regular epimorphisms and $Mono \; \mathcal{C}$ denotes the class of monomorphisms of $\mathcal{C}$. If $me$ is the $(RegEpi \; \mathcal{C},Mono \; \mathcal{C})$-factorization of a morphism $f$, then the pair $(Im\; f, m)$ (or, for short, $Im\; f$), where $Im\; f$ denotes the domain of $m$, is called the image of $f$. 

Before continue, recall that the category of points $Pt(\mathcal{C})$ of a category $\mathcal{C}$ is defined as follows: objects are pairs of morphisms $p:C\rightarrow C'$ and $s:C'\rightarrow C$ with $ps=1_{C'}$; morphisms $(p,s)\rightarrow (q,t)$ are pairs  of morphisms $u:C\rightarrow D$ and $v:C' \rightarrow D'$ from  $\mathcal{C}$ satisfying the conditions $qu=vp$ and $tv=us$; the composition of morphisms is defined componentwise. The codomain functor $\mathbf{F}: Pt(\mathcal{C})\rightarrow \mathcal{C}$ is a fibration \cite[Theorem 2.1.15]{BB}. 

A category $\mathcal{C}$ is called \textit{protomodular} if it has pullbacks of split epimorphisms, and all the inverse image functors of the fibration $\mathbf{F}$ reflect isomorphisms; the notion was introduced by Bourn \cite{B}. If a category is pointed and has pullbacks of split epimorphisms, then protomodulairity is equivalent to the split short five lemma, which states that, for a commutative diagram
\begin{equation}
\xymatrix{A\ar@{ >->}[r]^{m}\ar[d]_{u}&B\ar@{->>}[r]^{p}\ar[d]^{v}&C\ar[d]^{w}\\
A'\ar@{ >->}[r]^{m'}&B'\ar@{->>}[r]^{p'}&C'}
\end{equation}
with split epimorphisms $p$ and $p'$, and their kernels $m$ and $m'$ respectively, one has: if $u$ and $w$ are isomorphisms, then so is $v$ \cite[Proposition 3.1.2]{BB}.  In the case where $\mathcal{C}$ is a variety of universal algebras, there is a syntactical characterization of protomodularity due to Bourn and Janelidze \cite{BJ}.

A category is said to be \textit{homological} if it is pointed, regular and protomodular \cite{BB}. %Homological categories satisfy a number of nice properties of abelian categories. In particular, the classical homological lemmas (the short five lemma, the nine lemma, the snake lemma), as well as the Noether isomorphism theorems hold.
We refer the reader to the book \cite{BB} for various interesting properties of homological categories. Here, we provide only some of them, specifically, those that will be employed in Sections 3 and 4.
% One of them is the short five lemma, it is obtained from the split short five lemma, stated above, by replacing `split epimorphisms' in it by `regular epimorphisms'.

\begin{prop} \cite[Proposition 3.1.23, Corollary 4.1.3]{BB} \label{prop2.6}
Let $\mathcal{C}$ be a homological category. Then 
\vskip+1mm
(a) any normal monomorphism $m:C\rightarrowtail C'$ has a cokernel (denoted by $C'/C$), and is the kernel of that cokernel;
\vskip+1mm
(b) any regular epimorphism is the cokernel of its kernel.
\end{prop}

%\begin {theo} (First Noether isomorphism theorem) \cite[Theorem 4.3.10]{BB}
%Let $\mathcal{C}$ be a homological category. Consider two normal subobjects $H\subseteq C$ and $K\subseteq C$ with $H\subseteq K$. Then $H$ is a normal subobject of $K$, and $K/H$ is a normal subobject of $C/H$.  Moreover, the isomorphism $(C/H)/(K/H)\simeq C/K$ exists.
%\end{theo}

%Recall also the following properties of homological categories that we will need below. 

\begin{prop} \cite[Lemma 4.2.5(1, 2)]{BB} \label{prop2.7}
 For the commutative diagram 
\begin{equation}
\xymatrix{A\ar@{ >->}[r]^{m}\ar[d]_{u}&B\ar@{->>}[r]^{p}\ar[d]^{v}&C\ar[d]^{w}\\
A'\ar@{ >->}[r]^{m'}&B'\ar[r]^{p'}&C'}
\end{equation}
with a regular epimorphism $p$, and the kernels $m$ and $m'$ of $p$ and $p'$ respectively, one has: 

(a) if the left-hand square is a pullback, then the morphism $w$ is a monomorphism;

(b) if $u$ is an isomorphism, then the right-hand square is a pullback. 

(c) if $u$ and $w$ are regular epimorphisms, then so is $v$.
\end{prop}

\section{Suprema of normal subobjects}

As is well-known, suprema of finite families of subobjects  (in the partially ordered class of subobjects) of an object in an abelian category exist and have a simple construction. For an arbitrary homological category, the existence of the suprema is established in the case of finite families of normal subobjects. In the case of two normal subobjects with the zero intersection, the construction of supremum is given by the following theorem.

\begin{theo} (Couniversal property of products) \cite[Proposition 3.3.1, 3.3.2]{BB} \label{th2.5}
Let $\mathcal{C}$ be a homological category. For any objects $C$ and $C'$, the subobjects $l_C$ and $r_{C'}$ are normal, and their intersection is the zero subobject of $C\times C'$. Further, in the following diagram 
\begin{equation}
\xymatrix{C\ar@{ >->}[r]^{l_C}\ar@{ >->}[dr]_{m}&C\times C'\ar@{-->}[d]^{t}&C'\ar@{ >->}[l]_{r_{C'}}\ar@{ >->}[ld]^{m'}\\
&D&}
\end{equation}
if $m$ and $m'$ are normal subobjects whose intersection  is the zero subobject, then there exists a unique morphism $t$ making the diagram commutative. Moreover, $t$ itself is a normal subobject, which is the supremum of $m$ and $m'$ in the partially ordered class of subobjects of $D$.
\end{theo}

 Further, we provide some extracts from the proof of \cite[Theorem 4.3.12]{BB} in the form of the following two propositions (taking into account \cite[Lemma 4.3.1]{BB}). The first of them gives the construction of the supremum of two arbitrary normal subobjects.
\begin{prop} \label{prop3.1}
Let $\mathcal{C}$ be a homological category. Let $C$ be an object of $\mathcal{C}$, and $i_1:H\rightarrowtail C$ and $i_2:K\rightarrowtail C$ be its normal subobjects. Suppose $\varphi:C\rightarrow C/(H\cap K)$ is the cokernel of $H\cap K\rightarrowtail C$, while $\varphi(H)$ and  $\varphi(K)$ are the images of the compositions $\varphi i_1$ and $\varphi i_2$ resp. Then 
\vskip+1mm

(a) the intersection $H\cap K$ is a normal subobject of $H$, $K$, and $C$;
\vskip+1mm
(b) the subobjects $\varphi(H)$ and  $\varphi(K)$ of $C/(H\cap K)$ are normal, and their intersection $\varphi(H)\cap \varphi(K)$ is the zero object;
\vskip+1mm
(c) the supremum $H+K$ of $i_1$ and $i_2$ in the partially ordered class of subobjects of $C$ exists and is given by the monomorphism $i$ in the following diagram
\begin{equation}
\xymatrix{H\ar@{->>}[rrr]^{e_1}\ar@{ >-->}[dr]^{i'_1}\ar@{>->}[dddr]_{i_1}&&&\varphi (H)\ar@{ >->}[dddl]^{m_1}\ar@{>->}[dl]^{l_{\varphi(H)}}\\
&H+K\ar@{ >->}[dd]^{i}\ar@{->>}[r]^{\varphi_{H,K}\;\;\;\;}&\varphi(H)\times \varphi(K)\ar@{ >-->}[dd]^{t}&\\
\\
&C\ar@{->>}[r]^{\varphi \;\;\;\;\;\;\;}&C/(H\cap K)}
\end{equation}
here:  $t$ is the morphism induced by $m_1$ and $m_2$ according to the couniversal property of a product (see Theorem \ref{th2.5} and the claim (b) of this proposition); the internal square is a pullback; $i'_1$ is the unique morphism making the entire diagram commutative. The monomorphism $i'_2:K\rightarrowtail H+K$ is constructed similarly.
\end{prop}

\begin{prop} \label{prop3.2}
Under the conditions of Proposition \ref{prop3.1}, the composition $p_H\varphi_{H,K}$ (resp. $p_K\varphi_{H,K}$) is the cokernel of $i'_2$ (resp. $i'_1$), where $p_H$ (resp. $p_K$) is the projection $\varphi(H)\times \varphi(K)\twoheadrightarrow \varphi(H)$ (resp. $\varphi(H)\times \varphi(K)\twoheadrightarrow \varphi(K))$. 
\end{prop}

We are now ready to give the following statement that will be employed in Section 4.

\begin{prop} \label{prop3.3}
Let $\mathcal{C}$ be a homological category, and $C$ be an object in $\mathcal{C}$. Let  $i_1:H\rightarrowtail C$ and $i_2:K\rightarrowtail C$ be its normal subobjects, and $\pi_1$ and $\pi_2$ be their cokernels. Suppose $i: H+K\rightarrowtail C$ is their supremum in the partially ordered class of subobjects of $\mathcal{C}$. Then the following conditions are equivalent:
\vskip+1mm
(i) the composition $\pi_2 i_1: H\rightarrow C/K$ is a regular epimorphism;
\vskip+1mm
(ii) the monomorphism $i$ is an isomorphism;
\vskip+1mm
(iii) the composition $\pi_1 i_2: K\rightarrow C/H$ is a regular epimorphism.
\end{prop}

\begin{proof}
It suffices to prove the equivalence (i)$\Leftrightarrow$(ii). To this end, consider the diagram
\begin{equation}
\xymatrix{H\ar@{->>}[rrrr]^{e_1}\ar@{ >->}[dr]^{i'_1}\ar@{ >->}[dddr]&&&&\varphi (H)\ar@{ >->}[dddl]\ar@{>->}[dl]^{l_{\varphi(H)}}\\
&H+K\ar@{ >->}[dd]\ar@{->>}[rr]^{\varphi_{H,K}}&&\varphi(H)\times \varphi(K)\ar@{ >->}[dd]\ar@<0.7ex>@{->>}[dr]^{p_K}\ar@<0.7ex>@{->>}[ur]^{p_H}&\\
K\ar@{>->}[ur]^{\;\;\;\;\;\;\;\;\;\;\;\;\;\;\;\;\;\;\; \;\;\;\;\;\;\;\;\;\;\;\;\;\;\;\;\;\;\;\;\;\;\;\;\;\;\;\;i'_2}\ar[rrrr]^{\;\;\;\;\;\;\;\;\;\;\;\;\;\;\;\;\;\;\;\;\;\;i\;\;\;\;\;\;\;\;\;\;\;\;\;\;\;\;\;e_2\;\;\;\;\;\;\;\;\;\;\;\;\;\;\;\;\;\;\;\;\;\;\;\;\;\;\;\;\;\;\;\;\;\;\;\;\;t\;\;\;\;\;\;\;\;\;\;\;\;\;\;\;\;\;\;\;\;\;\;\;\;\;\;\;\;\;\;\;\;\;}\ar@{>->}[dr]_{i_2}&&&&\varphi(K)\ar@{>->}[ul]^{r_{\varphi_{K}}}\ar@{>->}[dl]^{m_2}\\
&C\ar@{->>}[dr]_{\pi_2}\ar@{->>}[rr]^{\varphi}&&C/(H\cap K)\ar@{-->}[dl]^{\sigma}\\
&&C/K}
\end{equation}
where $\sigma$ is the unique morphism with $\sigma\varphi = \pi_2$. 

We have $\pi_2 i_1=\sigma m_1 e_1$. Since $e_1$ is a regular epimorphism,  $\pi_2 i_1$ is a regular epimorphism if and only if $\sigma m_1$ is such. 

Note that
$$l_{\varphi_{H}}p_{H}\varphi_{H,K}i'_1=l_{\varphi_{H}}p_{H}l_{\varphi_{H}}e_1=l_{\varphi_{H}}e_1.$$
Since $l_{\varphi_{H}}$ is a monomorphism, this implies that $$p_{H}\varphi_{H,K}i'_1=e_1.$$

Further, Proposition \ref{prop3.2} implies that there exists a morphism $\psi:\varphi (H)\rightarrow C/K$ such that 
\begin{equation} \label{3.3}
\psi  p_{H}p_H \varphi_{H,K}=\pi_2 i.
\end{equation}
We have
$$\psi e_1=\psi p_H\varphi_{H,K}i'_1=\pi_2 i i'_1=\pi_2 i_1.$$
Since $e_1$ is an epimorphism, this implies that $\psi=\sigma m_1$. Taking into account equality (\ref{3.3}), we obtain that
\begin{equation}
\sigma m_1 p_{H} \varphi_{H,K} =  \pi_2 i,
\end{equation} 
and, therefore, the diagram
\begin{equation}
\xymatrix{K\ar@{ >-}[r]^{i'_2}\ar[d]_{1_K}&H+K\ar@{->>}[r]^{\;\;\;\;\;p_H\varphi_{H,K}}\ar@{ >->}[d]^{i}&\varphi(H)\ar[d]^{\sigma m_1}\\
K\ar@{ >->}[r]^{i_2}&C\ar@{->>}[r]^{\pi_2}&C/K}
\end{equation}
is commutative.  It follows that if $i$ is an isomorphism, then $\sigma m_1$ is a regular epimorphism.

For the converse, first note that, by Proposition \ref{prop3.2}, $p_H \varphi_{H,K}$ is a cokernel of $i_2'$. Since $i_2'$ is a normal monomorphism, it is the kernel of $p_H \varphi_{H,K}$. Therefore, if $\sigma m_1$ is a regular epimorphusm, then, by Proposition \ref{prop2.7}(c), $i$ is a regular epimorphism, and hence is an isomorphism. This completes the proof.

\end{proof}

% a morphism is a monomorphism if and only if its kernel is trivial;

\section{Radicals and their factorization systems on a homological category}
Let $\mathcal{C}$ be a homological category. A \textit{radical} on $\mathcal{C}$ is defined as a pair $(\mathbf{R},\mu)$ (or, for short, $\mathbf{R}$), where $\mathbf{R}$ is an endofunctor of $\mathcal{C}$ and $\mu$ is a natural transformation from $\mathbf{R}$ to the identity endofunctor such that, for any object $C$ of $\mathcal{C}$, $\mu_C$ is a normal monomorphism and $\mathbf{R}(C/\mathbf{R}(C))=0$.
\vskip+1mm
It is easy to observe that the full subcategory $\mathcal{X}$ of objects $X$ with $\mathbf{R}(X)=0$ is reflective in $\mathcal{C}$. The reflector $\mathbf{I}$ is given by $\mathbf{I}(C)=C/\mathbf{R}(C)$, while $\eta_C$ is the cokernel of $\mu_C$. 

%Let $\mathbb{E}$ be the class of morphisms of $\mathcal{X}$ which are regular epimorphisms in $\mathcal{C}$, while $\mathbb{M}$ be the class of monomorphisms of $\mathcal{X}$. 

\begin{lem} \label{lem4.1}
Let $\mathcal{C}$ be a homological category. Then a morphism of $\mathcal{X}$ is a regular epimorphism in $\mathcal{X}$ if and only if it is a regular epimorphism in $\mathcal{C}$. Moreover, the category $\mathcal{X}$ is homological. 
\end{lem} 

\begin{proof}
First note, that  the subcategory  $\mathcal{X}$ is closed under subobjects. Indeed, if $m:C\rightarrowtail X$ is a monomorphism with $X$ in $\mathcal{X}$, then there is $h:C/\mathbf{R}(C)\rightarrow X$ such that $h\eta_C=m$. This implies that $\eta_C$ is an isomorphism.

Let now $\mathbb{E}$ be the class of morphisms of $\mathcal{X}$ which are regular epimorphisms in $\mathcal{C}$. Then, applying Proposition \ref{prop2.6}(b), we obtain that $\mathbb{E}\subseteq RegEpi \; \mathcal{X}$. Moreover, for any morphism $f:X\rightarrow X'$ in $\mathcal{X}$, and its $(RegEpi \;\mathcal{C}, Mono \;\mathcal{C})$-factorization $f=me$ in $\mathcal{C}$, the domain of  $m$ is an object of $\mathcal{X}$. This implies that the pair $(\mathbb{E},Mono \;\mathcal{X})$ is a factorization system on $\mathcal{X}$.  Since $RegEpi \;\mathcal{X}\subseteq (Mono \; \mathcal{X})^{\uparrow}$, $\mathbb{E}= RegEpi \; \mathcal{X}$. Now it suffices to observe that $\mathcal{X}$ is closed under pullbacks.
%Let $f:X\rightarrow X'$ be a morphism in $\mathcal{X}$, and $f=me$ be its $(RegEpi,Mono)$-factorization in $\mathcal{C}$. The composition $\eta_C e$ is a regular epimorphism since the category $\mathcal{X}$ is regular. 

%There is a morphism $g$ with $m=g\eta_C$, where $C=codom \; e$. Show that $I(g)$ is a monomorphism. To this end, consider the commutative diagram
%\begin{equation}
%\xymatrix{V(C)\ar[d]_{V(C)}\ar@{>->}[r]&C\ar@{->>}[r]^{\eta_C}\ar@{>->}[d]^{m}&C/V(C)\ar[d]^{g}\\
%0\ar[r]&X'\ar[r]^{=}&X'}
%\end{equation}
%$V(m)$ is obviously a monomorphism, and hence an isomorphism. Theorem 2.7(a) implies that $g$ is a monomorphism.
\end{proof}

\begin{prop} \label{prop4.2}
Let $\mathcal{C}$ be a homological category. The class $\mathbf{I}^{-1}(RegEpi \; \mathcal{X})$ is that of morphisms $f:C\rightarrow C'$ such that the composition $\eta_{C'} f$ is a regular epimorphism. If the image of $f$ is a normal monomorphism, then $f\in \mathcal{C}$ if and only if $Im\; f+\mathbf{R}(C')=C'$. In particular, the latter criterion holds for any morphism $f$, if $\mathcal{C}$ is abelian.
\end{prop}

\begin{proof}
The first claim follows from the commutativity of the diagram
\begin{equation}
\xymatrix{C\ar@{->>}[r]^{\eta_C\;\;\;\;\;\;}\ar@{  >->}[d]_{f}&C/\mathbf{R}(C)\ar[d]^{I(f)}\\
C'\ar@{ >->}[r]^{\eta_{C'} \;\;\;\;\;}&C'/\mathbf{R}(C')}
\end{equation}
\noindent and the fact that the pair $(RegEpi \; \mathcal{C}, Mono \; \mathcal{C})$ is a factorization system on $\mathcal{C}$. 

For the second claim, consider any morphism $f$ and its $(RegEpi \;\mathcal{C}$, $Mono \;\mathcal{C})$-factorization $f=me$. Obviously, $\mathbf{I}(e)$ is a regular epimorphism. Therefore, $\mathbf{I}(f)$ is a regular epimorphism if and only if so is $\mathbf{I}(m)$. Therefore, (since the images of $f$ and $m$ coincide) it suffices to prove the claim in the case where $f$ is a monomorphism. But this immediately follows from Proposition \ref{prop3.3}.
\end{proof}

\begin{prop} \label{prop4.3}
Let $\mathcal{C}$ be a homological category. A morphism $m:C\rightarrow C'$ lies in $(\mathbf{I}^{-1}(RegEpi \; \mathcal{C}))^{\downarrow}$ if and only if $m$ is a monomorphism and satisfies the following condition: if $D$ is a subobject of $C'$ such that the composition of the inclusion $C\cap D\rightarrowtail D$ with the morphism $\eta_D$ is a regular epimorphism, then $D\subseteq C$.
\end{prop}

\begin{proof}
It suffices to apply Lemma \ref{lem2.4} with $(\mathbb{E}',\mathbb{M}')=(RegEpi \; \mathcal{C},$ $ Mono \; \mathcal{C}$).
\end{proof}

\begin{theo} \label{th4.4}
Let $\mathcal{C}$ be a complete well-powered homological category, and $\mathbf{R}$ be a radical. The pair of morphism classes $$(\mathbf{I}^{-1}(RegEpi \; \mathcal{X}), (\textbf{H}(Mono \; \mathcal{X})^{\uparrow\downarrow})$$ is a factorization system on $\mathcal{C}$. The class $\mathbf{I}^{-1}(RegEpi \; \mathcal{X})$ is described in Proposition 4.2, while the class $(\mathbf{H}(Mono \; \mathcal{X}))^{\uparrow\downarrow}=(\mathbf{I}^{-1}(RegEpi \; \mathcal{X}))^{\downarrow}$ is described in Proposition 4.3. 
\end{theo}

\begin{proof}
The claim follows from Lemma \ref{lem4.1} and Theorem \ref{th2.3}. 
\end{proof}

\begin{exmp} \label{exmp4.5}
There are many examples of complete well-powered homological categories. One of them is provided by the category $Top(\mathcal{V})$ of topological $\mathcal{V}$-algebras, where $\mathcal{V}$ is a pointed protomodular variety of universal algebras \cite[Theorem 4.6.5]{BB}. In view of Proposition  \ref{prop4.2}, note that, in this category, regular epimorphisms are precisely open surjective continuous homomorphisms \cite[Proposition 21]{BC}, while normal subobjects are precisely subobjects which are topological embeddings and are normal in $\mathcal{V}$ (or, equivalently, are ideals in the sense of the paper \cite{GU} by Gumm and Ursini. The characterization of ideals in pointed protomodular vareities is provided by \cite[Theorem 10.3.5]{CEL}). 
\end{exmp}

\begin{rem}
In view of Example \ref{exmp4.5} and Corollary \ref{th4.4}, note that there is another class of factorization systems on the category of topological $\mathcal{V}$-algebras, that naturally arises. Applying Theorem \ref{theo2.2}  to the forgetful functor
 $$\mathbf{I}:Top(\mathcal{V})\rightarrow \mathcal{V}$$ which is topological (and hence a fibration), we can construct a factorization system $(\mathbb{E}, \mathbb{M})$ on the category $Top(\mathcal{V})$, provided that a factorization system $(\mathbb{E}', \mathbb{M}')$  on the category $\mathcal{V}$ is given. We have: $f\in \mathbb{E}$ if and only if $\mathbf{I}(f)\in \mathbb{E}'$; $f\in \mathbb{M}$ if and only if $f\in \mathbb{M}'$ and $f$ is a topological embedding. Note that the fact that $(\mathbb{E}, \mathbb{M})$ is a factorization system can be easily verified directly.
\end{rem}

We say that a radical  $\mathbf{R}$ is \textit{idempotent} if the morphism $\mathbf{R}(\mu_C)(=\mu_{\mathbf{R}(C)})$ is an isomorphism, for any object $C$. 
\vskip+1mm
The following proposition generalizes  \cite[Proposition 5.6]{JT}.
\begin{prop} \label{prop4.5}
Let $\mathcal{C}$ be a homological category, and $\mathbf{R}$ be a radical on it. The reflection $\mathbf{I}\dashv \mathbf{H}$ is semi-left-exact if and only if $\mathbf{R}$ is idempotent.
\end{prop}

\begin{proof}

For "If" part, consider a pullback 
\begin{equation}
\xymatrix{A\ar[r]^{g}\ar[d]&X\ar[d]^{f}\\
C\ar[r]^{\eta_C}&I(C)}
\end{equation}
 with $\mathbf{R}(X)=0$. There is a morphism $\theta: \mathbf{R}(C)\rightarrow A$ with $f'\theta=\mu_C$ and $g\theta=0$:
\begin{equation}
\xymatrix{\mathbf{R}(A)\ar[d]_{\mathbf{R}(f')=f''}\ar@{ >->}[r]^{\;\;\;\;\mu_A}&A\ar[r]^{g}\ar[d]^{f'}&X\ar[d]^{f}\\
\mathbf{R}(C)\ar@{ >->}[r]^{\;\;\;\;\;\mu_C}\ar@{-->}[ur]^{\theta}&C\ar@{->>}[r]^{\eta_C\;\;\;\;\;\;}&C/\mathbf{R}C)}
\end{equation}
%The outer square of the following diagram is obviously commutative:
%\begin{equation}
%\xymatrix{VV(C)\ar[r]^{V(\varepsilon_C)}\ar[d]_{VV(f')}&V(C)\ar[d]^{f''=V(f')}\\
%VV(A)\ar[r]^{V(\varepsilon_A)}\ar[ur]^{V(\theta)}&V(A)}
%\end{equation}
One obviously has $\mathbf{R}(\mu_C)=\mathbf{R}(f')\mathbf{R}(\theta)$. Since $\mathbf{R}(\mu_C)$ is an isomorphism, $f''$ is a split epimorphism. 

Observe that $\mathbf{I}(g)\eta_A=g$. This implies that $g\mu_A=0=g\theta f''$. Besides, $$f'\theta f''=\mu_C f''=f'\mu_A.$$ This implies that $\theta f''=\mu_A$. Therefore, $f''$ is a monomorphism, and hence an isomorphism.

To show that $\mathbf{I}(g)$ is an isomorphism, we will show that $g$ is a cokernel of $\mu_A$.  To this end, it suffices to show  that $\mu_A$ is a kernel of $g$ due to the fact that $g$ is a regular epimorphism and Proposition \ref{prop2.6}(b). Consider any morphism $h:B\rightarrow A$ with $gh=0$. Since $\mu_C$ is a normal monomorphism, it is the kernel of $\eta_C$ by Proposition \ref{prop2.6}(a). Therefore, there is  a morphism $\alpha:B\rightarrow \mathbf{R}(C)$ with $f'h=\mu_C\alpha$. We have $\mu_Af''^{-1}\alpha=h$ as it follows from the equalities
$$f'\mu_Af''^{-1}\alpha=f'\theta\alpha=\mu_C\alpha=f'h$$ and $$g\mu_A f''^{-1}\alpha=g\theta \alpha=0=gh.$$  
%Since $\mathcal{C}$ is a regular category, $g$ is a regular epimorphism, and hence, by Theorem 2.5(b),  a cokernel of its kernel $\mu_A$.   

%\begin{equation}
%\xymatrix{C\ar[rr]^{g}\ar[dd]_{f'}\ar[dr]^{\eta_C}&&X\ar[dd]^{f}\\
%&C/V(C)\ar[ur]^{\beta}\\
%A\ar@{->>}[rr]&&A/V(A)}
%\end{equation}
 %is an isomorphism.  

For "Only if" part, consider the commutative diagram
\begin{equation}
\xymatrix{\mathbf{R}(C)\ar[rr]\ar[dr]^{\eta_{\mathbf{R}(C)}}\ar@{ >->}[dd]_{\mu_C}&&0\ar[dd]\\
&\mathbf{R}(C)/\mathbf{R}^{2}(C)\ar[ru]^{v}&\\
C\ar@{->>}[rr]^{\eta_C}&&C/\mathbf{R}(C)}
\end{equation}
\noindent It is obvious that the outer square is a pullback. Since the reflection is semi-left-exact, $\upsilon$ is an isomorphism. Hence $\eta_{\mathbf{R}(C)}$ is the zero morphism. This implies that its kernel $\mu_{\mathbf{R}(C)}$ is an isomorphism. 
\end{proof}

%\begin{proof}
%Consider the commutative diagram
%\begin{equation}
%\xymatrix{V(C)\ar@{>->}[r]^{\varepsilon_{C}}\ar@{>->}[d]_{V(m)}&C\ar@{->>}[r]^{\eta_C}\ar@{>->}[d]^{m}%&C/V\ar[d]^{I(m)}\\
%V(C')\ar@{>->}[r]^{\varepsilon_{C'}}&C'\ar@{->>}[r]^{\eta_C'}&C'/V(C')}
%\end{equation}
%with monomorphic $m$. First show that $\varepsilon_C$ is the kernel of $\eta_C'm$. Indeed, if %$h:K\rightarrow C$ is a morphism with $\eta_C'mh=0$, then there is $\alpha:K\rightarrow V(C')$ such $%%\varepsilon_Ch=\alpha$. Since the left-hand square is a pullback, there is $\beta$ with $%\varepsilon_C\beta=h$.

%Factorize now $I(m)$ as $I(m)=m'e$ with a regular epimorphism $e$ and a monomorphism $m'$. The monomorphism $\varepsilon_C$ is the kernel of the composition $gi\eta_C$ as is easy to verify. Since this composition is a regular epimorphism, it is a cokernel of $\varepsilon_C$. Therefore, $e$ is an isomorphism, and hence $I(m)$ is a monomorphism.
%\end{proof}

\begin{prop} \label{prop4.6}
Let $\mathcal{C}$ be a homological category, and $\mathbf{R}$ be an idempotent radical on it. The class $(\mathbf{H}(Mono \; \mathcal{X}))^{\uparrow\downarrow}$ is precisely the class of monomorphisms $C\rightarrowtail C'$ such that $\mathbf{R}(C')\subseteq C$.
\end{prop}

\begin{proof}
We are going to apply Proposition \ref{prop2.2}. If $\mathbf{I}(m)$ is a monomorphism and the right-hand square in the diagram
\begin{equation}
\xymatrix{\mathbf{R}(C)\ar@{ >->}[r]^{\mu_C}\ar[d]_{\mathbf{R}(m)}&C\ar@{->>}[r]^{\eta_C}\ar[d]^{m}&C/\mathbf{R}(C)\ar[d]^{\mathbf{I}(m)}\\
\mathbf{R}(C')\ar@{ >->}[r]^{\mu_{C'}}&C'\ar@{->>}[r]^{\eta_{C'}}&C'/\mathbf{R}(C')}
\end{equation}
 is a pullback, then obviously $m$ is a monomorphism, and there is $\theta:\mathbf{R}(C')\rightarrow C$ with $m\theta=\mu_{C'}$. For the converse, assume that $m$ is a monomorphism and the latter condition is satisfied. Then $\mathbf{R}(m)$ is a monomorphism as it follows from the commutativity of the left-hand square of diagram (4.4). On the other hand, it is a split epimorphism since $\mathbf{R}(m)\mathbf{R}(\theta)=\mathbf{R}(\mu_C)$. Therefore, $\mathbf{R}(m)$ is an isomorphism. This obviously implies that  the left-hand side of diagram (4.4) is a pullback. Then, by Proposition \ref{prop2.7}(a), the morphism $\mathbf{I}(m)$ is a monomorphism.  Moreover, by Proposition \ref{prop2.7}(b), the right-hand square in diagram (4.4) is a pullback. 
\end{proof}

\cite[Example 4.2]{CHK} implies that the requirement that $\mathbf{R}$ is idempotent cannot be removed from Proposition \ref{prop4.6}.

Corollary \ref{th2.1}, Lemma \ref{lem4.1}, and Proposition \ref{prop4.5} imply the following statement.

\begin{theo} \label{th4.9}
Let $\mathcal{C}$ be a homological category, and $\mathbf{R}$ be an idempotent radical on it. Then the pair of morphism classes $$(\mathbf{I}^{-1}(RegEpi \; \mathcal{X}), (\mathbf{H}(Mono \; \mathcal{X})^{\uparrow\downarrow})$$ is a factorization system on $\mathcal{C}$. The class $\mathbf{I}^{-1}(RegEpi \; \mathcal{X})$ is described in Propositions \ref{prop4.2}, while the class $(\mathbf{H}(Mono \; \mathcal{X}))^{\uparrow\downarrow}$ is described in Proposition \ref{prop4.6} (and, equivalently, in Proposition \ref{prop4.3}). 
\end{theo}

%\begin{rem}
%Under the conditions of Theorem \ref{th4.9}, the validity of the condition (b) in the definition of a factorization system is obvious for the case where a morphism $f$ has a  normal image $Im\;f$. As to the condition (c), it can be proved without Corollary \ref{th2.1}, but employing Propositions \ref{prop2.2}, \ref{prop4.5}, and  \ref{prop4.6}. To this end, first note, that the class of monomorphisms $m:C\rightarrowtail C'$ with $\mathbf{R}(C')\subseteq C$ is closed under pullback, as is easy to verify. Therefore, it suffices to show the existence of a diagonal morphism $\delta$ in diagram (2.1), for the case where $B=D$ and $\beta=1_B$ \cite{CJKP}. Consider the commutative diagram
%\begin{equation}
%\xymatrix{A\ar[r]^{e}\ar[d]^{\alpha}&B\ar@{->>}[r]^{\eta_B\;\;\;\;\;\;\;\;\;}&B/\mathbf{R}(B)\\
%C\ar[r]\ar@{ >->}[ru]_{m}&C/\mathbf{R}(C)\ar@{ >->}[ru]_{\mathbf{I}(m)}}
%\end{equation}
%Since the composition $\eta_{B}e$ is a regular epimorphism, so is $\mathbf{I}(m)$. But $\mathbf{I}(m)$ is a monomorphism, hence it is an isomorphism. Since $m$ is a trivial covering, the parallelogram on the right-hand part of the diagram is a pullback. Therefore, $m$ is an isomorphism. 
%\end{rem}

\vskip+2mm
\textit{Author's address:
Dali Zangurashvili, A. Razmadze Mathematical Institute of Iv. Javakhishvili Tbilisi State University},
\textit{Aleksidze Str., Lane II, Tbilisi 0193, Georgia, e-mail: dali.zangurashvili@tsu.ge}

\end{document}